\documentclass[12pt]{article}
\usepackage[cp1251]{inputenc}
\usepackage[english]{babel}
\usepackage{amssymb,amsfonts,amsmath}
\usepackage{amscd}
\newtheorem{Def}{Definition}[section]

\newtheorem{Th}{Theorem}[section]

\newtheorem{Lem}{Lemma}[section]
\newtheorem{Prop}{Propositon}[section]
\newtheorem{Cor}{Corollary}[section]

\newenvironment{Proof}
{\par\noindent{\bf Proof.}}
{\hfill$\scriptstyle\blacksquare$}

\title{Homotopy properties of the space $OS_{f}(X)$}
\author{Kh.~Kh.~Kurbanov \footnote{Academy of the Armed Forces of Uzbekistan, \textbf{e-mail:} qhamid\_83@mail.ru},
A.~Ya.~Ishmetov \footnote{Tashkent Institute of Architecture and Civil Engineering, \textbf{e-mail:} ishmetov\_azadbek@mail.ru}}

\begin{document}

\maketitle
\thispagestyle{empty}

\begin{center}

\end{center}
\begin{abstract}
For a given compact Hausdorff space $X$, we construct
the space $OS_{f}(X)$ of normed, order-preserving, weakly additive,
positively homogeneous and semi-additive functionals (for brevity,
semi-additive functionals) and it is proved that the hyperspace $\exp\, X$
of the space $X$ is a deformation retract of the constructed space.
Further we show that the shapes of the spaces $OS_{f}(X)$ and $\exp\, X$ coincide.
We establish that if the space $\exp\, X$ is contractible, then the space $OS_{f}(X)$ is also
contractible.\\

2010 \textit{MSC:} 54C65, 52A30.

\end{abstract}

\tableofcontents

\section{Introduction }

As it is known, the classical Krein-Milman theorem states that a convex compact set lying in a locally convex space coincides with the closure of the convex hull of the set of its extreme points. In the proof of the main achievements of this paper, this remarkable theorem occupies a central place. Although this theorem was established for linear spaces, it has recently been used for wider spaces (see, for example, \cite{Davletov2009}, \cite{Zait2006}, \cite{Zait2019fujma}, \cite{Zait2020agt}, \cite{Radul1999}).

Consider a compact Hausdorff space $X$, a Banach algebra $C(X)$ of all continuous functions \ $\varphi\colon\, X \to \mathbb{R}$, provided with pointwise algebraic operations and a $\sup$-norm, that is, the norm $\|\varphi\|=\{|\varphi(X)|:\, x\in X\}$.
For each $c\in\mathbb{R}$, $c_X$ denotes the constant function defined by the formula $c_X(x)=c$, $x\in X$. Let $\varphi$, $\psi \in C(X)$. The inequality $\varphi\leq\psi$ means that $\varphi(x)\leq\psi(x)$  for all $x\in X$.

\begin{Def}\label{semiadditivefunc}
{\rm
A functional $\mu\colon\,C(X)\to\mathbb{R}$ is called:
\begin{enumerate}
\item \textit{weakly additive}, if, for all $c\in\mathbb{R}$ and $\varphi \in C(X)$ the equality $\mu(\varphi +c_X)=\mu(\varphi)+c$ holds;
\item \textit{order-preserving } if, for any pair of functions $\varphi$, $\psi \in C(X)$, the inequality $\varphi\leq\psi$ implies $\mu(\varphi)\leq\mu(\psi)$;
\item \textit{normed}, if $\mu(1_X)=1$;
\item \textit{positively homogeneous}, if $\mu(t\varphi) = t\mu(\varphi)$ for all $\varphi \in C(X)$, $t\in\mathbb{R}_{+}$, where $\mathbb{R}_{+}=[0,+\infty)$;
\item \textit{semiadditive}, if $\mu(\varphi + \psi) \leq \mu(\varphi) + \mu(\psi)$ for all  $\varphi$, $\psi \in C(X)$.
\end{enumerate}}
\end{Def}

For a compact Hausdorff space $X$, by $OS(X)$ we denote the set of all functionals $\nu\colon\, C(X)\to\mathbb{R}$ that satisfy the above five conditions, and for brevity, these functionals are called semiadditive functionals.

The set $OS(X)$ provided with the point-wise convergence topology. Note, the sets of the view
\begin{gather*}
\langle\mu;\, \varphi_1,\, \dots,\, \varphi_n;\, \varepsilon\rangle=\{\nu\in OS(X):\, |\mu(\varphi_i)-\nu(\varphi_i)<\varepsilon,\ i=1,\, \dots, \,n\},
\end{gather*}
form a base of neighbourhoods of a functional $\mu \in OS(X)$ in the point-wise convergence topology, where $\varphi_i \in C(X)$, $i=1,\,\dots,\, n$, $n\in \mathbb{N}$, $\varepsilon>0$.

It is clear, that if $\mu$, $\nu\in OS(X)$, then $\alpha\mu + \beta\nu \in OS(X)$, where $\alpha$, $\beta\geq0$, $\alpha + \beta = 1$.
Moreover, the following statement is true.

\begin{Th}{\rm
\cite{Davletov2009}. For any compact Hausdorff space $X$, the space $OS(X)$ is a convex compact with respect to the point-wise convergence topology.}
\end{Th}

Let $X$ and $Y$ be compact Hausdorff spaces, $f\colon\, X\to Y$ be a continuous map. By the formula
\begin{gather*}
OS(f)(\mu)(\varphi)=\mu(\varphi\circ f),\qquad \mu \in OS(X),
\end{gather*}
we define a map $OS(f)\colon\, OS(X)\to OS(Y)$, $\varphi \in C(Y)$.

The operation $OS$ defines a covariant functor which acts in the category $\mathfrak{Comp}$ of compact Hausdorff spaces and their continuous maps. Note, $OS$ is a normal functor. Obviously, for every compact Hausdorff space $X$, the space $P(X)$ of probability measures (i.~e. linear, non-negative, normed functionals) is a subspace of $OS(X)$.

Let $A$ be a nonempty subset of the space $P(X)$ of probability measures on $X$, $\varphi \in C(X)$. Then $|\mu (\varphi)|<\|\varphi\|$ for any $\mu \in A$, and therefore a number set $\{\mu (\varphi):\, \mu \in A\}$ is bounded from above.
Consequently, for every $\varphi \in C(X)$ there is a number
\begin{gather}\label{nu_A}
\nu_A(\varphi)=\sup \{\mu(\varphi):\, \mu \in A\}.
\end{gather}

For compact Hausdorff space $X$ by $\exp\, X$ we denote the hyperspace of $X$, that is, the space of all nonempty closed subsets of $X$ provided by the Vietoris topology (see \cite{ZaitJuma2019ujma}, \cite{ZaitJuma2020} for more details). For each $F\in \exp\, X$ we define a functional $\mu_{F}\colon\, C(X)\to \mathbb{R}$ as following
\begin{gather}\label{mu_{F}}
\mu_{F}(\varphi)=\max_{x\in F} \varphi(x),\ \ \varphi \in C(X).
\end{gather}
Clearly, $\mu_{F}$ is a weakly additive, order-preserving, normed, positively-homogenous and semiadditive functional. The correspondence $F\mapsto \mu_{F}$  is one-to-one. Therefore, we can the set $F$ identify to the functional $\mu_{F}$. Thus, $\exp\, X\subset OS(X)$.

Let $K$ be a closed subset of some locally convex space $E$. By $cc(K)$ we denote a set consisting of all convex closed subsets of $K$ and on $cc(K)$ consider the topology induced from $\exp\, K$.

\begin{Th}\label{tauhomeo}
{\rm \cite{Davletov2009}.
Let $X$ be a compact Hausdorff space. Then spaces $OS(X)$ and $cc(P(X))$ are homeomorphic. This homeomorphism $\tau\colon\, cc(P(X))\to OS(X)$ may by define by the rule
\begin{gather*}
\tau (A)=\nu_A, \qquad\, A\in cc(P(X)).
\end{gather*}
}
\end{Th}

From theorem \ref{tauhomeo} and the above mentioned Krein-Milman theorem we get the following statement.

\begin{Cor}\label{Cortauhomeo}
{\rm
Let $X$ be a compact Hausdorff space. Then spaces  $OS(X)$ and $cc(P(X))$ are homeomorphic. The homeomorphism $\tau\colon\, cc(P(X))\to OS(X)$ may be defined as 
\begin{gather*}
\tau (A)=\nu_{\operatorname{ext}\,A}, \qquad A\in cc(P(X)),
\end{gather*}
where $\operatorname{ext}\,A$ is a set of all extreme points of a convex compact set $A \subset P(X)$.
}
\end{Cor}

By reformulating the definition 18 from \cite{Shchepin1981}, we introduce the concept of a support of semiadditive functional.
A \textit{support} of~$\mu\in OS(X)$ is a closed subset~$\mbox{supp}\,\mu\subset X$ such that relations  ~$A\supset\mbox{supp}\,\mu$ and~$\mu\in OS(A)$ are equivalent for each closed $A\subset X$. For a functor the~$OS$ the support exists for every~$\mu\in OS(X)$ and it defines as
\[
\mbox{supp}\,\mu = \cap \left\{ A\subset X:\overline{A}=A,\,\,\mu \in OS(A) \right\},
\]
here~$\overline{A}$ is the closure of~$A$.

For a compact Hausdorff space $X$ and a positive integer $n$ by $OS_n (X)$ we denote a set of all functionals $\mu \in OS(X)$ for which $|\operatorname{supp}\,\mu | \leq n$. $OS_n (X)$ consider as a subspace of the space $OS(X)$. Put
\begin{gather*}
OS_\omega(X)=\underset{i=1}{\overset{\infty}{\cup}} OS_{n}(X).
\end{gather*}

A functional $\mu \in OS_{\omega} (X)$ is called as semiadditive functional with the finite support. Theorem \ref{tauhomeo} and the Corollary \ref{Cortauhomeo} imply that for each semiadditive functional $\mu$ with the finite support there exists only unique convex closed set of $A\in cc(P(X))$ that
\begin{gather*}
\mu = \nu_{A} = \nu_{\operatorname{ext}\,A}.
\end{gather*}
It is clear that
\begin{gather*}\label{suppOS=cupsuppP}
\operatorname{supp}\, \mu = \underset{\xi\in A}{\cup}\operatorname{supp}\,\xi.
\end{gather*}

At the same time, if each element $\xi\in A$ is a probability measure with the finite support, assume $\operatorname{supp}\, \xi = \{x_{\xi1},\,\dots,\, x_{\xi\,n_{\xi}}\}$, then
\begin{gather*}
\xi = \underset{i=1}{\overset{n_{\xi}}{\Sigma}}\alpha_{\xi\,i}\delta_{x_{\xi\,i}},
\end{gather*}
where $\underset{i=1}{\overset{n_{\xi}}{\Sigma}}\alpha_{\xi\,i} = 1$, $\alpha_{\xi\,i}\geq 0$, $i = 1,\, \dots,\, n_{\xi}$. Hence, the formula (\ref{nu_A}) can be written as
\begin{gather*}
\mu(\varphi)=\sup \left\{\underset{i=1}{\overset{n_{\xi}}{\Sigma}}\alpha_{\xi\,i}\delta_{x_{\xi\,i}}(\varphi):\, \xi \in A\right\}. \tag{\ref{nu_A}$^{\, \prime}$}
\end{gather*}

The following set was introduced by E.~V.~Shchepin.
\begin{multline*}
P_{f}(X)=\Big\{\mu \in P_\omega (X):\, \mbox{ if }\,  \mu = \underset{i=1}{\overset{n}{\Sigma}}\alpha_i \delta_{x_i}, \\
\mbox{ then there exists }\,  i_0 \in \{1,\dots,\, n\} \mbox{ such that }\,   \alpha_{i_0}\geq 1-\frac{1}{n+1}\Big\}.
\end{multline*}
For a compact Hausdorff space $X$, we define the following set
\begin{gather*}\label{OS_{f}}
OS_{f}(X)=\big\{\nu_{A} \in OS(X):\, \operatorname{ext}\,A\subset P_{f}(X)\big\}.
\end{gather*}

\begin{Lem}\label{expsubsetOS}
{\rm
For an arbitrary compact Hausdorff space $X$ and every closed $F\subset X$, we have a $\mu_{F} \in OS_{f}(X)$.
}
\end{Lem}
\begin{Proof}
Let $F\in \exp\, X$. Let us consider the functional $\mu_{F}$ defined by equality (\ref{mu_{F}}). Then $\mu_{F} = \nu_{\{\delta_{x}:\, x\in F\}}$. Hence $\mu_{F} \in OS_{f}(X)$.

\end{Proof}

From the proved Lemma \ref{expsubsetOS}, in particular it follows that $OS_{f}(X) \neq \varnothing$.

\begin{Prop}\label{OSiscomp}
{\rm
The topological space $OS_{f} (X)$, equipped with the point-wise convergence topology, is a compact Hausdorff space.
}
\end{Prop}
\begin{Proof}
From the homeomorphism $OS_{f}(X) \cong cc(P_{f}(X))$ immediately follows the proof of the proposition.

\end{Proof}

Lemma \ref{expsubsetOS} and Proposition \ref{OSiscomp} give the following statement.
\begin{Cor}
{\rm
For an arbitrary compact Hausdorff space $X$, its hyperspace $\exp\, X$ is homeomorphic to some closed subset of the compact Hausdorff space $OS_{f}(X)$.
}
\end{Cor}

Remind the following concept. Let $\mathfrak{C} = \{\mathfrak{O},\, \mathfrak{M}\}$ and $\mathfrak{C}^{\, \prime}={\{\mathfrak{O}^{\, \prime},\, \mathfrak{M}^{\, \prime}}\}$ be two categories, where $\mathfrak{O}$, $\mathfrak{O}^{\, \prime}$ are the classes of objects, and $\mathfrak{M}$, $\mathfrak{M}^{\, \prime}$ are the classes of morphisms. A map $F\colon\, \mathfrak{C}\to \mathfrak{C}^{\, \prime}$, transforming objects to objects, and morphisms to morphisms, is said to be a \textit{covariant functor} acting from the category $\mathfrak{C}$ into the category $\mathfrak{C}^{\, \prime}$ if it satisfies the following conditions:
\begin{itemize}
\item[$F1)$] For every morphism $f\colon\, X\to Y$ from the category $\mathfrak{C}$, the morphism $F(f)$ acts from $F(X)$ to $F(Y)$;
\item[$F2)$] $F(\operatorname{id}_{X}) = \operatorname{id}_{F(X)}$ for all $X\in\mathfrak{O}$;
\item[$F3)$] $F(f\circ g)=F(f)\circ F(g)$  for every pair of morphisms $f$ and $g$ from $\mathfrak{M}$.
\end{itemize}

\begin{Prop}\label{OSfisfunctor}
{\rm
The construction $OS_{f}$ forms a covariant functor acting in the category $\mathfrak{Comp}$ of compact Hausdorff spaces and their continuous maps.
}
\end{Prop}
\begin{Proof}
At first, for a given map $f\colon\, X\to Y$ we show that $OS(f)(OS_{f}(X))\subset OS_{f}(Y)$.
Let $\mu \in OS_{f}(X)$. Then there exists an $A\subset P(X)$ such that $\operatorname{ext}\, A\subset P_{f}(X)$ and $\mu = \nu_{\operatorname{ext}\,A}$. But, then from the results of work \cite{Zait2019} it follows that $P_{f}(f)(\operatorname{ext}\,A) = \operatorname{ext}\, f(A)\subset P_{f}(Y)$. Therefore, $OS_{f}(\mu) = OS_{f}(\nu_{\operatorname{ext}\,A}) = \nu_{\operatorname{ext}\, f(A)}\in OS_{f}(Y)$. Now let us define a map $OS_{f}(f)\colon\, OS_{f}(X)\to OS_{f}(Y)$ as the restriction $OS_{f}(f) = OS(f)|_{OS_{f}(X)}$. Since $OS(f)$ is continuous \cite{Davletov2009}, its restriction $OS_{f}(f)$ is also continuous. Thus, $OS_{f}$ satisfies the condition $F1)$.

Let $\operatorname{id}_X\colon\, X\to X$ be the identity map. For every $\mu\in OS_{f}(X)$ we have
\begin{gather*}
OS_{f}(\operatorname{id}_X)(\mu)(\varphi )=\mu (\varphi\circ \operatorname{id}_X)=\mu(\varphi), \qquad \varphi\in C(X).
\end{gather*}
Since $\mu$ and $\varphi$ are arbitrary, then it becomes $OS_{f}(\operatorname{id}_X)(\mu)=\mu$ (the condition $F2)$ carried out).

Let us show that $OS_{f}$ preserves the map composition. Let $X$, $Y$, $Z$ be compact Hausdorff spaces and $f\colon\, X\to Y$, $g\colon\, Y\to Z$ be continuous maps. For $\mu \in OS_{f} (X)$ and $\varphi \in C(Z)$ we have
\begin{multline*}
OS_{f} (g\circ f)(\mu)(\varphi)=\mu (\varphi \circ (g\circ f))=\mu ((\varphi \circ g)\circ f) =OS_{f} (f)(\mu)(\varphi \circ g)=\\
= OS_{f}(g)\circ OS_{f}(f)(\mu)(\varphi),
\end{multline*}
i.~e.~$OS_{f}(g\circ f)= OS_{f}(g)\circ OS_{f}(f)$ (the condition $F3)$ is established).

\end{Proof}

Thus, the construction $OS_{f}$ transforming compact Hausdorff spaces to compact Hausdorff spaces, and continuous maps of compact Hausdorff spaces in continuous maps of compact Hausdorff spaces, forms a covariant functor acting in the category of compact Hausdorff spaces and their continuous mappings.

In the works of \cite{Shchepin1981} E.~V.~Schenpin, \cite{Fedorchuk1991} V.~V.~Fedorchuk, \cite{Zait2019} A.~A.~Zaitov, \cite{Juraev1989} T.~F.~Jurayev were studied the functor $P_{f}$ a traditional probability analogue of the built functor $OS_{f}$. The work of \cite{ZaitIshm2019} A.~A.~Zaitov and A.~Ya.~Ishmetov was devoted to the idempotent analogue $I_{f}$.

In the present work, we establish that the functor $OS_{f}$ is a normal in the category of compact Hausdorff spaces and their continuous maps. Further we prove that the hyperspace $\exp\, X$ of the compact Hausdorff space $X$ is a deformation retract of the space $OS_{f}(X)$. We also show that the shapes of $OS_{f}(X)$ and $\exp\, X$ are the same. We prove that if $\exp\, X$ is contractible, then $OS_{f}(X)$ is also contractible. Obtained results strictly differ from the above-described authors' results, because in previous works have been established relationships between the compact Hausdorff space $X$ and $I_{f}(X)$ or $P_{f}(X)$.

But all these three functors have one remarkable property: the degrees of the functors $OS_{f}$, $P_{f}$ and $I_{f}$ are infinite.
But this phenomenon immediately generates another difference: the functors $P_{f}$ and $I_{f}$ are with a finite support, and the construction $OS_{f}$ is a functor with an infinite support.

Recall the concept of the degree of functors. Let~$X$ be a compact Hausdorff space, $F$ is a functor and~$x\in F(X)$. A \textit{degree of the point}~$x$ is a smallest positive integer $n$ such that~$x$ belongs to~$F(f)F(K)$ for some map~$f\colon K\to X$ of $n$-point space in~$K$ (\cite{Shchepin1981}, Definition 16). If does not exist such finite $n$ , then degree of~$x$ is considered as \textit{infinite}. \textit{Degree of a functor}  of~$F$ is the maximum of degrees of various points~$x\in F(X)$ for various compact Hausdorff spaces~$X$ and it denotes by~$\mbox{deg} F$.
\\

\section{Normality of the functor $OS_{f}$}

A functor~$F\colon \mathfrak{Comp}\to \mathfrak{Comp}$ acting in the category of compact Hausdorff spaces and their continuous maps is said to be \textit{normal} if it satisfies the following conditions (\cite{Shchepin1981}, Definition 14):
\begin{enumerate}
\item $F$ is continuous ($F(\lim S) = \lim F(S)$),
\item $F$ preserves weight ($wX = wF(X)$),
\item $F$ is monomorphic (i.~e. preserves the injectivity of maps),
\item $F$ is epimorphic (i.~e. preserves the surjectivity of maps),
\item $F$ preserves the intersections ($F(\bigcap\limits_{\alpha} X_\alpha) = \bigcap\limits_{\alpha} (F(X_\alpha)$),
\item $F$ preserves the preimages ($F(f^{-1}) = F(f)^{-1}$),
\item $F$ preserves a point and an empty set ($F(\mathbf{1}) = \mathbf{1}$, $F(\varnothing) = \varnothing$).
\end{enumerate}

Let us decipher this definition. Let~$S=\{X_\alpha,\, p_\alpha^\beta;\, \mathfrak{A}\}$ be the inverse spectrum of compact Hausdorff spaces, $\lim S = \lim\limits_{\leftarrow} S$ is its limit. According to Kurosh theorem, the inverse spectrum limit of non-empty compact Hausdorff spaces is non-empty (\cite{Fedorchuk2006}, Theorem 3.13) and is a compact Hausdorff space (\cite{Fedorchuk2006}, Proposition 3.12). Under the impact of the functor~$F$ to the compact Hausdorff spaces~$X_\alpha$ and to the maps~$p_\alpha^\beta$, $\alpha,\,\beta\in\mathfrak{A}$, $\alpha\prec\beta$, the inverse spectrum~$F(S)=\{F(X_\alpha),\, F(p_\alpha^\beta);\, \mathfrak{A}\}$ is formed. Let~$\lim F(S)$ be the limit of this spectrum. Condition 1 requires that the equality~$F(\lim S) = \lim F(S)$ has to hold. For a topological space~$X$ by~$wX$ we denote its weight, that is, the smallest of the base powers of the space~$X$. Condition 2 requires that the weights of the compact Hausdorff spaces~$X$ and~$F(X)$ be equal. The monomorphism of the functor~$F$ (condition 3) allows us to consider~$F(A)$ as the subspace~$F(X)$ for the closed~$A \subset X$. The identity of~$F(A)$ with the subspace of~$F(X)$ is done by the idendity~$F(i_A)$, where~$i_A\colon A \to X$ is the identity. Condition 4 requires that if~$f\colon X\to Y$ is a continuous map ``onto'', then~$F(f)\colon F(X)\to F(Y)$  is also a continuous map ``onto''. For a monomorphic functor~$F$, conditions 5 and 6 are deciphered as follows: for any family~$\{X_\alpha\}$ of closed subsets of an arbitrary compact Hausdorff space~$X$, the equality~$F(\bigcap\limits_{\alpha} X_\alpha) = \bigcap\limits_{\alpha}F(X_\alpha)$ has to hold (condition 5); For every continuous map of~$f\colon X\to Y$ and every closed~$B$ in~$Y$, the equality~$F(f^{-1}(B)) = F(f)^{-1}F(B)$ (condition 6) is true. The point preservation condition means that~$F$ takes a one-point space to a one-point space.

\begin{Prop}
{\rm
The functor $OS_{f}$ preserves the weight of infinite compact Hausdorff spaces,~i.~e. for every infinite compact Hausdorff space the equality $w(OS_{f}(X))=w(X)$ holds.
}
\end{Prop}
\begin{Proof}
From the relations $X \cong \delta (X) \subset OS_{f}(X) \subset OS(X)$ and the equality $w(OS(X))=w(X)$ established in \cite{Davletov2009} follows the required equality.

\end{Proof}
\begin{Prop}
{\rm
$OS_{f}$ is a monomorphic functor,~i.~e., it preserves the injectivity of mappings of compact Hausdorff spaces.
}
\end{Prop}
\begin{Proof}
Let $\mu_1$, $\mu_2 \in OS_{f} (X)$, $\mu_1\neq \mu_2$. Owing to the injectivity of the map $f$ there exists a function $\varphi \in C(Y)$, such that $\mu_1 (\varphi \circ f)\neq \mu_2(\varphi \circ f)$. Hence $OS_{f} (f)(\mu_1)(\varphi) = \mu_1 (\varphi \circ f)\neq \mu_2(\varphi \circ f)=OS_{f} (f)(\mu_2)(\varphi)$.

\end{Proof}
\begin{Prop}
{\rm
If $f\colon\, X\to Y$ is a continuous map ``onto'', then $OS_{f} (f)\colon\,OS_{f} (X)\to OS_{f} (Y)$ is also a continuous ``onto'' map.
}
\end{Prop}
\begin{Proof}
The continuity of the map $OS_{f} (f)$ had shown in Proposition \ref{OSfisfunctor}. Since the for a surjective map $f\colon\, X\to Y$ the map $OS(f)$ is surjective \cite{Davletov2009}, then its restriction $OS_{f} (f)$ is also surjective.

\end{Proof}
\begin{Prop}
{\rm
Functor $OS_{f}\colon\, \mathfrak{Comp}\to \mathfrak{Comp}$ preserves
\begin{itemize}
\item[$a)$] a point,	
\item[$b)$] the empty set.
\end{itemize}
}
\end{Prop}
\begin{Proof} a) Let $x\in X$. By definition, we have $OS_{f}(\{x\}) = \{\delta_x \}$.

b) Let $X=\varnothing$. Then $C(X)=\varnothing$. Consequently, $\mathbb{R}^{C(X)} =\mathbb{R}^\varnothing=\varnothing$. From here we get $OS_{f} (\varnothing)\subset \varnothing$.

\end{Proof}
\begin{Prop}
{\rm
If $A$ is a closed subset of a compact Hausdorff space $X$, then $OS_{f} (A) \subset OS_{f} (X)$.
}
\end{Prop}
\begin{Proof}
Let $A$ be closed in $X$ and $\mu \in OS_{f} (A)$. Then the functional $\mu$ is concentrated on $A$. Owing to the definition of the concept of the support it is equivalent to $\operatorname{supp}\, \mu\subset A$. Then $\operatorname{supp}\,\mu\subset X$, from where $\mu \in OS_{f} (X)$.

\end{Proof}
\begin{Prop}
{\rm
If $f\colon\, X\to Y$ is a continuous map between compact Hausdorff spaces and $B\subset Y$, then $OS_{f} (f^{-1} (B)) = OS_{f} (f)^{-1} (OS_{f} (B))$.
}
\end{Prop}
\begin{Proof}
Let $\mu \in OS_{f} (f^{-1} (B))$. By the definition, this means that $\mu\in OS_{f} (X)$ and $\operatorname{supp}\,\mu \subset f^{-1} (B)$. Consequently, $f(\operatorname{supp}\,\mu)\subset B$. Therefore $\operatorname{supp}\, OS_{f} (f)(\mu)\subset B$. From here $OS_{f} (f)(\mu)\in OS_{f} (B)$, i.~e. $ \mu \in OS_{f} (f)^{-1} (OS_{f} (B))$.

Inversely, let $\mu \in OS_{f} (f)^{-1} (OS_{f} (B))$. Then $OS_{f} (f)(\mu)\in OS_{f} (B)$,~i.~e. $\operatorname{supp}\,OS_{f} (f)(\mu)\subset B$. Consequently,  $f(\operatorname{supp}\,\mu)\subset B$. This means that $\operatorname{supp}\,\mu \subset f^{-1} (B)$, whence $\mu \in OS_{f} (f^{-1} (B))$.

\end{Proof}

Let $\{X_\alpha,\, p_\alpha^\beta;\,  A\}$ be an inverse spectrum indexed by the elements of the set $A$ and consisting of compact Hausdorff spaces. By $\lim X_\alpha$ we denote the limit of this spectrum, and by $p_\alpha\colon\, \lim X_\alpha \to X_\alpha, \alpha \in A$, the limit projections. Inverse spectrum $\{X_\alpha, \,p_\alpha^\beta;\,  A\}$ generates the inverse spectrum $\{OS_{f} (X_\alpha),\, OS_{f} (p_\alpha^\beta);\,  A\}$, which limit is denoted by $\lim OS_{f}(X_\alpha)$, and the limit projections by $pr_\alpha\colon\, \lim OS_{f} (X_\alpha) \to OS_{f} (X_\alpha)$. The maps $OS_{f} (p_\alpha)\colon\, OS_{f} (\lim X_\alpha)\to OS_{f} (X_\alpha),\colon\, \alpha \in A$, generate the map $R_{OS_{f}}\colon\, OS_{f} (\lim X_\alpha)\to \lim OS_{f} (X_\alpha)$.

\begin{Prop}
{\rm
The functor $OS_{f}$ is continuous,~i.~e., the map $R_{OS_{f}}\colon\, OS_{f} (\lim X_\alpha)\to \lim OS_{f} (X_\alpha)$ is a homeomorphism. }
\end{Prop}
\begin{Proof}
Since taking an affine combination and taking a closure are continuous operations, it follows from the continuity \cite{Davletov2009} of the functor $OS$ that $R_{OS_{f}}\colon\, OS_{f} (\lim X_\alpha)\to \lim OS_{f} (X_\alpha)$ is a homeomorphism.

\end{Proof}
\begin{Prop}
{\rm
The functor $OS_{f}$ preserves the intersection, i.~e., for any pair of closed subsets $A$, $B$ of a compact Hausdorff space $X$, we have
\begin{gather*}
OS_{f} (A\cap B)=OS_{f} (A)\cap OS_{f} (B).
\end{gather*}
}
\end{Prop}
\begin{Proof}
The inclusion $OS_{f} (A\cap B)\subset OS_{f} (A)\cap OS_{f} (B)$ is clear. If $\mu \in OS_{f} (A)\cap OS_{f} (B)$, then $\operatorname{supp}\,\mu \subset A$ and $\operatorname{supp}\,\mu \subset B$; consequently, $\operatorname{supp}\,\mu \subset A\cap B$. From here $\mu\in OS_{f} (A\cap B)$, i.~e. $OS_{f} (A\cap B)\supset OS_{f} (A)\cap OS_{f} (B)$.

\end{Proof}

Thus, the following main result of the section is proved.

\begin{Th}
{\rm
$OS_{f}\colon\, \mathfrak{Comp}\to \mathfrak{Comp}$ is a normal functor.
}
\end{Th}

\section{The contractibility of the space of semiadditive functionals}

In this section, we will establish that if for a given Hausdorff compact space~$X$ its hyperspace~$\exp\,X$ is a contractible compact, then~$OS_f(X)$ is also a contractible compact.

A subset~$Y$ of a topological space~$X$  is (\cite{Borsuk1976}, p. 14) a \textit{retract} of~$X$ if there exists a map~$r\colon X\to Y$ (called a \textit{retraction} of~$X$ into~$Y$) such that the restriction~$r|_Y\colon Y\to Y$ is the identity map ~$id_Y\colon Y\to Y$ (i.~e.,~$r(y)=y$ for all~$y\in Y$).
If ~$r\colon X\to Y$ is a retraction and there exists a homotopy~$h\colon X\times[0,\,1]\to Y$ such that ~$h(x,0)=x$, $h(x,1)=r(x)$, for all~$x\in X$, then~$r$ is a \textit{deformation retraction}, and~$Y$ is a \textit{deformation retract} of the space~$X$.
A deformation retraction~$r\colon X\to  F$ is a \textit{strongly deformation retraction} if, for the homotopy ~$h\colon X\times[0,\,1]\to  Y$, we have~$h(x,t)=x$ for all~$x\in F$ and all~$t\in[0,\,1]$.
A space~$Y$ is an \textit{absolute retract} (and they write~$Y\in AR$) if, for every homeomorphism~$h$ that maps~$Y$ onto a closed subset~$hY$ of any space~$X$, the set~$hY$ is a retract of the space~$X$. A space~$Y$ is called an \textit{absolute neighborhood retract} (and they write~$Y\in ANR$) if, for every homeomorphism~$h$ mapping~$Y$ onto a closed subset~$hY$ of any space~$X$, there exists a neighborhood~$U$ of the set~$hY$ (in~$X$) such that~$hY$ is a retract for~$U$.

Let~$X$ and~$Y$ be two compact sets lying in the metrizable spaces~$M$ and~$N$, respectively, where~$M,\,N\in AR$.
The sequence of mappings~$f_{k}\colon M\to N$, $k=1,\,2,\,\dots,$  is called (\cite{Borsuk1976}, p. 17) the \textit{fundamental sequence} from~$X$ to~$Y$, if for each neighborhood~$V$ of the compactum~$Y$ (in~$N$) there exists a neighborhood~$U$ of the compactum~$X$ (in~$M$) such that
\[
f_{k}|_{U}\simeq f_{k+1}|_{U} \quad \mbox{in} \quad V \quad \mbox{for almost all} \quad k=1,\,2,\,\dots.
\]
Here, ``for almost all'' means ``for all but finite number''. The relation~$f_{k}|_{U}\simeq f_{k+1}|_{U}$ means that there exists a homotopy~$\varphi_{k}\colon U\times[0,\,1]\to  V$ such that~$\varphi_{k}(x,\,0)=f_{k}(x)$ and~$\varphi_{k}(x,\,1)=f_{k+1}(x)$ for all~$x\in U$. This fundamental sequence is denoted by ~$\left\{f_{k},\,X,\,Y\right\}_{M,\,N}$ or shortly by~$\mathbf{f}$, and they write~$\mathbf{f}\colon X\to Y$ in~$M$, $N$. They say that the fundamental sequence~$\mathbf{f}=\left\{f_{k},\,X,\,Y\right\}_{M,\,N}$ is generated by the map~$f\colon X\to Y$ if~$f_{k}(x)=f(x)$ for all~$x\in X$ and for all~$k=1,\,2,\,\dots$.

Let~$X$ and~$Y$ be closed subsets of metrizable~$AR$-spaces~$M$ and~$N$, respectively.
They say (\cite{Borsuk1976}, p. 29) that the spaces~$X$ and~$Y$ are \textit{fundamental equivalent} (with respect to~$M$, $N$) if there exist two fundamental sequences~$\mathbf{f}\colon X\to Y$ and~$\mathbf{g}\colon Y\to X$ such that ~$\mathbf{g}\mathbf{f}=id_{X,\, M}$ and~$\mathbf{f}\mathbf{g}=id_{Y,\, N}$. The fundamental equivalence relation is equivalence relation. Therefore, the class of all spaces splits into pairwise disjoint classes of spaces, which are called \textit{shapes} (\cite{Borsuk1976}, p. 31). Consequently, two spaces belong to the same shape if and only if they are fundamental equivalent. The shape containing the space~$X$ is called the \textit{shape of the space}~$X$  and it denotes by~$Sh(X)$. The concept of shape is topological, i.~e., two homeomorphic spaces have the same shape. It is known that for two absolute neighborhood retracts~$A$ and~$B$, the equality~$Sh(A)=Sh(B)$ holds if and only if they are homotopy equivalent.

Take an arbitrary functional~$\nu_{A} \in OS_{f}(X)$. Note that each probability measure $\xi$ with finite support, say, $\{x_{1}, \, \dots,\, x_{n}\}$, is represented in the form of an affine combination $\xi = \overset{n}{\underset{i=1}{\Sigma}}\alpha_{i}\delta_{x_{i}}$ of Dirac measures $\delta_{x_{i}}$, $i = 1,\, \dots,\, n$, uniquely ( see, for example, \cite{Fedorchuk1991}).

Let $\operatorname{ext}\,A = \{\xi_{s}:\, s\in S\}$. By definition, $\operatorname{ext}\,A\subset P_{f}(X)$. It is clear that the support of each $\xi_{s} \in \operatorname{ext}\,A$ is finite; suppose $|\operatorname{supp}\, \xi_{s}|= n_{s}$. Let $\xi_{s} = \overset{n_{s}}{\underset{i=1}{\Sigma}}\alpha_{s,i}\delta_{x_{s,i}}$, $s\in S$. By construction, for each $\xi_{s}$ there exists $i(s)\in \{1, \, \dots,\, n_{s}\}$ such that $\alpha_{k,i(s)} \geq 1-\frac{1}{n_{s} + 1}$. To the functional $\nu_{A}$ we associate the set $F_{A}: = \{x_{1,i(s)}:\, s\in S\}$. The defined correspondence~$OS_{f}(X)\to \exp\, X$ is denoted by~$r_{e}^{O} = r_{\exp\, X}^{OS_{f}(X)}$. The map~$r_{e}^{O}\colon OS_{f}(X)\to  \exp\, X$ is defined correctly. The construction of the map~$r_{e}^{O}$, easily implies that ~$r_{e}^{O}(\mu_{F}) = F$ for each~$F\in\exp\, X$, i.~e., under the map~$r_{e}^{O}$, the points of the space~$\exp\, X$ are fixed points, where $\mu_{F}$ is the functional defined by equality (\ref{mu_{F}}). Consequently,~$r_{e}^{O}$ is a retraction, and the set $\exp\, X$ is a retract of the set $OS_{f}(X)$.

We will establish a stronger statement. To prove it, we identify the set $F\in \exp\, X$ with the functional $\mu_{F} \in OS_{f}(X)$.

\begin{Th}\label{expisdeforret}
{\rm
For an arbitrary Hausdorff compact space $X$, the set $\exp\, X$ is a strongly deformation retract of the Hausdorff compact space $OS_{f} (X)$.
}
\end{Th}
\begin{Proof}
Consider a map ~$h\colon OS_{f}(X)\times[0,\,1]\to  OS_{f}(X)$ defined by the formula
\begin{gather*}
h(\mu,\,t)=h_{t}(\mu)=(1 - t)\cdot\mu + t\cdot r_{e}^{O}(\mu), \qquad (\mu,\,t)\in OS_{f}(X)\times[0,\,1].
\end{gather*}

It is easy to verify that the map~$h$ is well defined. Moreover, $h_{0}=\operatorname{id}_{OS_{f}(X)}$ and~$h_{1}=r_{e}^{O}$, i.~e.,~$h$ is the homotopy connecting maps~$\operatorname{id}_{OS_{f}(X)}$ and~$r_{e}^{O}$.
Further, we have
\[
h(\mu_{F},\, t) = (1 -t )\cdot\mu_{F} + t\cdot r_{e}^{O}(\mu_{F})=\mu_{F},
\]
i.~e.,~$h_{t}(\mu_{F})=\mu_{F}$ for all~$F\in\exp\, X$ and~$t\in[0,\,1]$. Thus, $\exp\, X$ is a strongly deformation retract of the compact ~$OS_{f}(X)$.

\end{Proof}

From Theorem \ref{expisdeforret} and statement (5.4) from \cite{Borsuk1976} (p. 32) we obtain

\begin{Cor}
{\rm
For an arbitrary compact Hausdorff space~$X$, we have
\[
Sh(\exp\,X) = Sh(OS_f(X)).
\]
}
\end{Cor}

Recall (\cite{Borsuk1971}, p. 29) that a~$A\subset X$ is \textit{contractible in the space}~$X$ to \textit{to the set}~$B\subset X$ if the embedding~$i_{A}\colon A\to  X$ is homotopic to some map~$f\colon A\to  X$ such that ~$f(A)\subset B$. If in this case~$B$ consists of only one point, then they say that~$A$ is \textit{contractible in}~$X$.

Clearly, if there exists a homotopy~$h\colon A\times [0;\,1]\to  A$, such that~$h(y,0)=i_{A}$, and~$h(y,1)=\{\text {point}\}$, then~$A$ is contractible in ~$X$.

A space~$X$ is called (\cite{Borsuk1971}, p. 31) \textit{locally contractible} at a point~$x_{0}\in X$ if every neighborhood~$U$ of the point~$x_{0}$ contains a neighborhood~$U_{0}$ contractible in~$U$ to a point. A space~$X$ is called \textit{locally contractible} if it is locally contractible at each of its points.

\begin{Th}
{\rm
If for a Hausdorff compact space $X$ its hyperspace $\exp\, X$ is contractible, then the space $OS_{f} (X)$ is also contractible.
}
\end{Th}
\begin{Proof}
We show more: the functor~$OS_{f}$ preserves the homotopy of maps. Let~$h_{0}$, $h_{1}\colon X\to  Y$ be homotopical maps, $h\colon X\times[0,\,1]\to  Y$ be the homotopy connecting the maps~$h_{0}, \ \ h_{1}$, i.~e.~$h(x,\,0)=h_{0}(x)$, $h(x,\,1)=h_{1}(x)$. The embedding~$i_{t_{0}}\colon X\times\{t_{0}\}\to  X\times I$ defined by the equality~$i_{t_{0}}(x,\,t_{0})=(x,\,t_{0})$,  $x\in X$, defines the embedding~$OS_{f}(i_{t_{0}})\colon OS_{f}(X\times\{t_{0}\})\to OS_{f}(X\times I)$. But, for every~$t_{0}\in[0,\,1]$, the space ~$OS_{f}(X\times\{t_{0}\})$ is naturally homeomorphic to~$OS_{f}(X)\times\{t_{0}\}$. This homeomorphism can be realized, as it is easy to see, using equality (\ref{nu_A}$^{\, \prime}$) and the correspondence~$\mu_{t_{0}}\leftrightarrow(\mu,\,t_{0})$, where for each $\{\xi_{s}:\, s\in S\}\subset P_{f}(X)$:
\begin{align*}
\mu(\varphi)&=\sup \left\{\underset{i=1}{\overset{n_{s}}{\Sigma}}\alpha_{s\,i}\delta_{x_{s\,i}}(\varphi):\, s \in S\right\}, \qquad \varphi\in C(X),\\
\mu_{t_{0}}(\phi)&=\sup \left\{\underset{i=1}{\overset{n_{s}}{\Sigma}}\alpha_{s\,i}\delta_{(x_{s\,i},\, t_{0})}(\phi):\, s \in S\right\}, \qquad \phi\in C(X\times \{t_{0}\}).
\end{align*}

We now define a map~$OS_{f}(h)\colon OS_{f}(X)\times[0,\,1]\to  OS_{f}(Y)$ by the equality
\[
OS_{f}(h)\left(\sup \left\{\underset{i=1}{\overset{n_{s}}{\Sigma}}\alpha_{s\,i}\delta_{x_{s\,i}}:\, s \in S\right\},\,t\right) = \sup \left\{\underset{i=1}{\overset{n_{s}}{\Sigma}}\alpha_{s\,i}\delta_{h(x_{s\,i},\,t)}:\, s \in S\right\}.
\]
We have
\begin{gather*}
OS_{f}(h)\left(\sup \left\{\underset{i=1}{\overset{n_{s}}{\Sigma}}\alpha_{s\,i}\delta_{x_{s\,i}}:\, s \in S\right\},\,0\right) = \sup \left\{\underset{i=1}{\overset{n_{s}}{\Sigma}}\alpha_{s\,i}\delta_{h(x_{s\,i},\,0)}:\, s \in S\right\} =\\
= \sup \left\{\underset{i=1}{\overset{n_{s}}{\Sigma}}\alpha_{s\,i}\delta_{h_{0}(x_{s\,i})}:\, s \in S\right\} = OS_{f}(h_{0})\left(\sup \left\{\underset{i=1}{\overset{n_{s}}{\Sigma}}\alpha_{s\,i}\delta_{x_{s\,i}}:\, s \in S\right\}\right),
\end{gather*}
\begin{gather*}
OS_{f}(h)\left(\sup \left\{\underset{i=1}{\overset{n_{s}}{\Sigma}}\alpha_{s\,i}\delta_{x_{s\,i}}:\, s \in S\right\},\,1\right) = \sup \left\{\underset{i=1}{\overset{n_{s}}{\Sigma}}\alpha_{s\,i}\delta_{h(x_{s\,i},\,1)}:\, s \in S\right\} =\\
= \sup \left\{\underset{i=1}{\overset{n_{s}}{\Sigma}}\alpha_{s\,i}\delta_{h_{1}(x_{s\,i})}:\, s \in S\right\} = OS_{f}(h_{1})\left(\sup \left\{\underset{i=1}{\overset{n_{s}}{\Sigma}}\alpha_{s\,i}\delta_{x_{s\,i}}:\, s \in S\right\}\right),
\end{gather*}
i.~e.~$OS_{f}(h)(\mu,\,0) = OS_{f}(h_{0})(\mu)$  and~$OS_{f}(h)(\mu,\,1)=OS_{f}(h_{1})(\mu)$ for each~$\mu\in OS_{f}(X)$. In other words, $OS_{f}(h)$ is the homotopy connecting ~$OS_{f}(h_{0})$  and~$OS_{f}(h_{1})$ maps. Thus, the functor~$OS_{f}$ preserves the homotopy of maps.

\end{Proof}

\end{document}